# Solvable Matrix Groups and Burnside's Problem
## S. Bachmuth

## 1. Introduction

In sections 2 & 3 we describe four matrix groups that we designate $F(\mathcal{R})$, $F(S)$, $F(\mathcal{R}[t, t^{-1}])$ and $F(S[t, t^{-1}])$. These 2x2 matrix groups are the subjects of the four part Theorem 1. The solvable group $F(S[t, t^{-1}])$ is the subject of Theorems 2 & 3 in section 4. Section 5 has applications to the 2-generator Burnside groups $B(q)$ of prime power exponent q. A consequence of Theorem 4 is the solvability of $B(q)$.

Theorem 1:

(**i**) There exists a (rank 2) free group $F(\mathcal{R}[t, t^{-1}])$ of 2x2 matrices with entries in the Laurent polynomial ring $\mathcal{R}[t, t^{-1}]$, where $\mathcal{R}$ itself is the Laurent polynomial ring with integer coefficients $\mathcal{R} = Z[x, x^{-1}, y, y^{-1}]$.

(**ii**) The homomorphic image $F(\mathcal{R})$ of $F(\mathcal{R}[t, t^{-1}])$, obtained by setting $t = 1$ in each element of $F(\mathcal{R}[t, t^{-1}])$, is a group of matrices with entries in $\mathcal{R}$ that is isomorphic to $F/F''$, the free metabelian group of rank 2.

(**iii**) For each positive integer n, there exists a quotient ring $S = S(n)$ of $\mathcal{R}$ such that if n is a prime-power, $F(S)$ is isomorphic to $F/F''F^n$, the free metabelian Burnside group of exponent n and rank 2.

(**iv**) The group $F(S[t, t^{-1}])$, obtained by replacing the ring $\mathcal{R}$ by the ring $S$ in $F(\mathcal{R}[t, t^{-1}])$, is solvable if (and only if) the integer n is a prime power. (The solvability class approaches infinity as the prime power n does.)

We thus have the following commutative square where the horizontal maps $\alpha$ and $\underline{\alpha}$ come from the ring homomorphism $\mathcal{R} \to S$ and the vertical maps are those determined by sending t to the identity in each matrix.

$$\begin{array}{ccc} F(\mathcal{R}[t, t^{-1}]) & \xrightarrow{\alpha} & F(S[t, t^{-1}]) \\ \downarrow & & \downarrow \\ F(\mathcal{R}) & \xrightarrow{\underline{\alpha}} & F(S) \end{array}$$

Note that for prime-power exponents, $\underline{\alpha}: F(\mathcal{R}) \to F(S)$ is the Burnside map on the free metabelian group. It is likely true for arbitrary exponents, but our proof is valid only for prime-power exponents. (Map means group homomorphism.)

The group $G = F(S[t, t^{-1}])$ in part (iv) of Theorem 1 is the subject of section 4. This matrix group depends on the ring $S = S(n)$ and we may write $G = G(n)$ if we wish to emphasize this dependence on n. When n = q is a prime power, $G = G(q) = F(S[t, t^{-1}])$ is solvable and the results in [3] enable one to determine the exact solvability class. As will become apparent, the class is determined by the ideal in $\mathcal{R}$ that is used to construct $S$, an ideal which depends on q. The solvability class of G(q) increases with increasing q. The question of why G(q) should be solvable for any q is particularly interesting since G(q) has many elements of infinite order - including one of it's two generators! Theorem 2 is the outgrowth of this question and provides our answer.

G(q) is a group of matrices over the ring $S[t, t^{-1}]$ and eliminating t in these matrices (i.e., factoring by the ideal t -1) yields the group F($S$), a group isomorphic to the rank 2 metabelian Burnside group of exponent q. In other words, the largest possible two generator metabelian group of exponent q. The motivation for constructing $S$ as a quotient of the (Laurent polynomial) ring $\mathcal{R}$ was to produce F($S$). At this point it is natural to construct the ring $S[t, t^{-1}]$ from the ring $\mathcal{R}[t, t^{-1}]$ to produce the group $G(q) = F(S[t, t^{-1}])$ as a quotient of the free group $F(\mathcal{R}[t, t^{-1}])$. Thus a matrix in G(q) has exponent q modulo G″ and we would expect a reasonable number of elements of G(q) to have exponent q. (For example, an element without t in any entry clearly has exponent q.) But much more is true. Leaning heavily on Theorem B in [3], it transpires that an element of G(q) has either infinite order or order dividing q, and those with exponent q are more ubiquitous than expected. For example, the commutator subgroup of G(q) has exponent q. But even here the picture is incomplete. Theorem 2 (in section 4) is a precise statement as to which elements have exponent q.

In this way, the groups G(q) offer strong evidence for the existence of higher commutator identities in groups of exponent q, a result well known for small q. How else to explain the solvability of G(q)? The ideal in $\mathcal{R}$ used to construct $S$ was chosen to produce elements of exponent q. Yet by means of Theorem B in [3], very high commutator identities unfold for large q. Sections 2 and 4 have the details.

The proofs of (**i**), (**ii**) and (**iii**) of Theorem 1 are in section 2,, the proof of (**iv**) is in section 3. Section 4 contains a detailed discussion of the groups $G(q) = F(S[t, t^{-1}])$. In section 5 we show that the Burnside map $\beta: F(\mathcal{R}[t, t^{-1}]) \to B(q) \cong F/F^q$ coincides with the map $\alpha: F(\mathcal{R}[t, t^{-1}]) \to F(S[t, t^{-1}])$ on the elements of $\mathcal{R}$. Specifically, we show

<u>Theorem</u> 4: Suppose g is an element of $F(\mathcal{R}[t, t^{-1}])$. Then in the entries of g, $\alpha$ and $\beta$ are the same on $\mathcal{R}$., i.e., both send the elements of $\mathcal{R}$ into their respective elements of $S$.

Theorem 4 should be understood in conjunction with Theorem 1 which tells us that the image of $\alpha$ is solvable. The proof of Theorem 4 is self contained and elementary.

Comments:

1) If a group G is such that $G/G'' \cong F/F''F^n$, and n is not a prime power, then $G/G''$ is not nilpotent. Most important for us is that the n-cyclotomic ideal introduced in section 2 is easier to deal with when n is a prime-power.

2) Although $F(S[t, t^{-1}])$ is solvable if (and only if) $n = q = p^e$ is a prime power, it is not nilpotent. Since this fact is not required anywhere, we will omit the proof which we find much harder than merely showing that $F(S[t, t^{-1}])$ is solvable. That $F(S[t, t^{-1}])$ has elements of infinite order is easy to show.

3) Dark and Newell [4] have shown that the exact nilpotency class of $F/F''F^q$ is $e(p^e - p^{e-1})$ when F has rank 2 and $q = p^e$ is a prime power. Their result can be used to give an alternative method for calculating the solvability class of $F(S[t, t^{-1}])$. For our purposes, it is more straightforward to use Theorem B in Bachmuth, Heilbronn, and Mochizuki [3].

4) The computation of solvability bounds for $F(S[t, t^{-1}])$ essentially reduces to the knowledge of the nilpotency class of $F/F''F^q$. In considering these metabelian groups of prime power exponent $q = p^e$, one can draw upon Theorem B in Bachmuth, Heilbronn, and Mochizuki [3] when the rank $r \leq p+1$; and Theorem 3 in Gupta, Newman, and Tobin [6] when $e > 1$ and $r > (p+2)(e-1)$. The aforementioned theorem of Dark and Newell [4] can also be used when $r = 2$.

In this paper we restrict ourselves to two generators. The extension to more than two generators requires nothing new conceptually and will be taken up at a later time. However, complications may arise in generalizing Lemma 3 (iii) of section 2 since the Jacobi identity now comes into play (see Theorem 3 of [1] ).

## 2. Metabelian matrix groups

Let F be the free group of rank 2 generated by the following matrices $M_1$ and $M_2T$ with entries in the Laurent polynomial ring $Z[x, x^{-1}, y, y^{-1}, t, t^{-1}]$.

$$M_1 = \begin{vmatrix} 1 & 1-y \\ 0 & x \end{vmatrix} \qquad M_2T = \begin{vmatrix} yt & 0 \\ 1-xt & 1 \end{vmatrix}$$

That F is free can be seen by taking $x = 1$ and $y = t = -1$ and appealing to a well known theorem. cf. Sanov [8].

Notation: We set $\mathcal{R} = Z[x, x^{-1}, y, y^{-1}]$ and write $F = F(\mathcal{R}[t, t^{-1}])$. The units of $\mathcal{R}$ are the monomials $x^i y^j$ and $-x^i y^j$ where $i, j$ are integers. We call $x^i y^j$ the positive units.

Lemma 1: (i) The image of F upon setting $t = 1$ is isomorphic to the free metabelian group $F/F''$ generated by the following matrices $M_1$ and $M_2$ with entries in $\mathcal{R} = Z[x, x^{-1}, y, y^{-1}]$. We put $F(\mathcal{R}) = gp\langle M_1, M_2 \rangle$.

$$M_1 = \begin{vmatrix} 1 & 1-y \\ 0 & x \end{vmatrix} \qquad M_2 = \begin{vmatrix} y & 0 \\ 1-x & 1 \end{vmatrix}$$

(ii) Each element M of $F(\mathcal{R})$ has the form $M = uI + N$, where I is the identity matrix, u is a (positive) unit in $\mathcal{R}$, and the 2x2 matrix N has the form

$$N = \begin{vmatrix} \lambda_1(1-x) & \lambda_1(1-y) \\ \lambda_2(1-x) & \lambda_2(1-y) \end{vmatrix}, \text{ where } \lambda_1, \lambda_2 \text{ in } \mathcal{R} \text{ satisfy}$$

$\lambda_1(1-x) + \lambda_2(1-y) = 1-u$. We will omit the I and write $M = u + N$.

Proof: (i) $F(\mathcal{R})$ is isomorphic to the group of inner automorphisms of the rank 2 free metabelian group as described in [1]. ([1] provides a description of the group of all IA-automorphisms of the free metabelian group which contains the inner automorphisms.)

(ii) A proof is contained in [2], section 2, Proposition(ii). Because of the importance of this Lemma, we shall reproduce the proof in an appendix to this paper.

Notation: Let $v = (1-x, 1-y)$ and write N as above in the form $N = [\lambda_i v]$, where $\lambda_i v$ denotes the $i^{th}$ row of N.

The lower central series for a group G is defined inductively by $G_1 = G$ and for $j > 1$ the $j^{th}$ term is $G_j = [G, G_{j-1}]$.

**Lemma 2:** If $M = u + N$ is in $F(\mathcal{R})$, where $N = [\lambda_i v]$, then

    (i) $M^n = u^n + (1 + u + u^2 + \cdots + u^{n-1})N$ for any integer $n > 1$.

    (ii) If $M$ is in the commutator subgroup of $F(\mathcal{R})$, then $\lambda_1$ and $\lambda_2$ are in $\Sigma$, the augmentation ideal of $\mathcal{R}$. More generally, if $M$ is in the $j^{th}$ term of the lower central series of $F(\mathcal{R})$, then $\lambda_1$ and $\lambda_2$ are in $\Sigma^{j-1}$.

**Proof:** Note that $N^2 = (1-u)N$ since $\lambda_1(1-x) + \lambda_2(1-y) = (1-u)$.

(i) From Lemma 1(ii) and the fact that $N^2 = (1-u)N$, we have
$M^2 = (u + N)^2 = u^2 + 2uN + (1-u)N = u^2 + (1+u)N$. Now proceed by induction.
$M^n = M^{n-1}M = \{u^{n-1} + (1+u+\cdots+u^{n-2})N\}(u + N) = u^n + (1+u+u^2+\cdots+u^{n-1})N$ again using $N^2 = (1-u)N$.

(ii) Choose a basis for $G = F(\mathcal{R})$ modulo $G_{j+1}$ and, using induction on $j$, compute using this basis. A convenient basis (H. Neumann [6], chapter 3, section 6) is the basic commutators that are not in the second derived group. Notice that if $M$ is in the commutator subgroup, then $u = 1$ and $\lambda_1(1-x) + \lambda_2(1-y) = 0$. Hence $\lambda_1$ and $\lambda_2$ are in the augmentation ideal of $\mathcal{R}$. In fact, $[M_2, M_1] = I + [\lambda_i v]$, where $\lambda_1 = -(1-y)$ and $\lambda_2 = (1-x)$. This is the first case and the start of the induction. It is straightforward to check that if the basic commutator $[M_2, M_1, M_2, \ldots, M_2, M_1, \ldots, M_1] = I + [\lambda_i v]$, where $M_1$ and $M_2$ occur a total of $a + 1$ and $b + 1$ times respectively, then $\lambda_1 = -(1-y)^{b+1}(1-x)^a$ and $\lambda_2 = (1-y)^b(1-x)^{a+1}$. Furthermore, conjugation of a commutator by $M = u + N$ becomes multiplication (of the non identity term) by $u$. Thus, Lemma 2(ii) follows.

    Before proceeding we need to introduce more notation. $q = p^e$ will always denote a prime power and we let $\phi(q) = p^e - p^{e-1}$, the Euler $\phi$-function. Let $\mathcal{I}(n)$ be the q-cyclotomic ideal, the ideal in $\mathcal{R}$ generated by all elements $1 + u + u^2 + \cdots + u^{n-1}$ where $u$ is a positive unit in $\mathcal{R}$. We denote by $\mathcal{S}$ the quotient ring. $\mathcal{R}/\mathcal{I}(n)\Sigma$. In Lemma 3(i)&(iii) we require that $n$ be a prime-power.

**Lemma 3:** (i) $\Sigma^{e\phi(q)} \subseteq \mathcal{I}(q)$ in $\mathcal{R}$. Equality holds if (and only if) $e=1$, that is, $q = p^e$ is prime. Furthermore $\Sigma^{e\phi(q)-1} \not\subseteq \mathcal{I}(q)$.

    If $e \geq 2$, $0 \leq j \leq e-1$ and $k \geq e(p^e - p^{e-1}) - j(p^{e-1} - p^{e-2})$, then $p^j \Sigma^k \subseteq \mathcal{I}(q)$.

    (ii) The group generated by the matrices $M_1$ and $M_2$ over $\mathcal{R}/\Sigma^c$ is isomorphic to $F/F''F_c$ the free metabelian nilpotent group of class $c$.

    (iii) The group $F(\mathcal{S})$ generated by the matrices $M_1$ and $M_2$ over $\mathcal{R}/\mathcal{I}(q)\Sigma$ is isomorphic to $F/F''F^q$, the free metabelian Burnside group of prime power exponent $q$.

Proof: (i) This is part of Theorem B in [3] and the Proposition on page 241 of [3]. Theorem B is the case j = 0 and will suffice for our purposes. However, for computations when e >1, the extra information in [3] is desirable. (Arbitrary finite rank is needed for a later paper. For now we are concerned only with rank r = 2.)

(ii) This follows from Lemmas 1 and 2(ii). For clearly a matrix of $F(\mathcal{R})$ which is in $F(\mathcal{R})_c$ has entries congruent to I modulo $\Sigma^c$ by Lemma 2(ii). Conversely, suppose $M = u + [\lambda_i v] \equiv I \mod \Sigma^c$. Then $u = 1$, $\lambda_i v$ are in $\Sigma^c$ and $\lambda_1(1-x) + \lambda_2(1-y) = 0$. Thus $\lambda_1 = \lambda(1-y)$, $\lambda_2 = -\lambda(1-x)$ and $\lambda$ is in $\Sigma^{c-2}$. From the proof of Lemma 2(ii), it is clear that M is in $F(\mathcal{R})_c$.

(iii) This follows from Lemmas 1 and 2(i). For again it is clear by Lemma 2(i) that a matrix in $F(\mathcal{R})'$ which is a $q^{th}$ power has entries congruent to I modulo $\mathcal{I}(q)\Sigma$. Conversely, suppose $M = u + [\lambda_i v] \equiv I$ modulo $\mathcal{I}(q)\Sigma$. Assume first that M is in the commutator subgroup. Then $u = 1$, $\lambda_i$ are in $\mathcal{I}(q)$, and $\lambda_1(1-x) + \lambda_2(1-y) = 0$. Thus the $\lambda_i$ are in the augmentation ideal and again $\lambda_1 = \lambda(1-y)$, $\lambda_2 = -\lambda(1-x)$ and both $\lambda(1-y)$ and $\lambda(1-x)$ are in $\mathcal{I}(q)$. By Theorem B of [3] or the Theorem of Dark and Newell [4], we conclude that $\lambda(1-y)$ and $\lambda(1-x)$ are in $\Sigma^{e\phi(q)}$ and thus M is in $F(\mathcal{R})^q$.

Any element M in $F(\mathcal{R})$ can be written $M = M_1^i M_2^j C$, where C is in the commutator subgroup of $F(\mathcal{R})$. Then $M = x^i y^j + [\lambda_i v]$ for suitable $\lambda_i$. By the hypothesis, i and j must each be a multiple of q. Hence, by Lemma 2(i), $M_1^i$ and $M_2^j$ are each congruent to the identity. So we may assume M is in the commutator subgroup. This completes the proof of Lemma 3.

## 3. Power series representations of F

We return to our free group $F = F(\mathcal{R}[t, t^{-1}]) = gp\langle M_1, M_2T \rangle$ defined at the beginning of the previous section. To deal with solvable images of $F(\mathcal{R}[t, t^{-1}])$, we shall describe $F(\mathcal{R}[t, t^{-1}])$ as a power series in $(t-1)^k$, $k = 0,1,2,...$, where the coefficients are 2x2 matrices over $\mathcal{R} = Z[x, x^{-1}, y, y^{-1}]$. Namely, let

$$T = \begin{vmatrix} t & 0 \\ 1-t & 1 \end{vmatrix} \quad \text{and} \quad S = \begin{vmatrix} 1 & 0 \\ -1 & 0 \end{vmatrix}$$

Then,
$M_2T = M_2(I + (t-1)S) = M_2 + (t-1)M_2S$

$(M_2T)^{-1} = (I + (t-1)S)^{-1}M_2^{-1} = (I - (t-1)S + (t-1)^2S^2 - (t-1)^3S^3 + \cdots)M_2^{-1}$

$= M_2^{-1} - (t-1)SM_2^{-1} + (t-1)^2SM_2^{-1} - (t-1)^3SM_2^{-1} + \cdots$

Since $F(\mathcal{R}[t, t^{-1}]) = gp\langle M_1, M_2T \rangle$, we can state formally

Lemma 4: Any element f in $F(\mathcal{R}[t, t^{-1}])$ has a representation of the form

(*)  $f = M_f + (t-1)A_1 + (t-1)^2 A_2 + (t-1)^3 A_3 + \cdots$

where $M_f \in F(\mathcal{R}) = gp\langle M_1, M_2 \rangle$ and the $A_i$ as well as $M_f$ are matrices over $\mathcal{R}$.

We know from Lemma 1(i) that $F/F''$ is the image of F via the map $f \to M_f$. The next couple of lemmas state further consequences of this description of $F(\mathcal{R}[t, t^{-1}])$.

Notation: Recall that $\Sigma$ denotes the augmentation ideal of $\mathcal{R}$. If all the entries of a matrix A are in $\Sigma^i$ we say that A is in $\Sigma^i$.

Lemma 5: Suppose $f \in F$ as in (*) is in F' the commutator subgroup of F. Then all the $A_i$ (the coefficients of $(t-1)^i$) are in $\Sigma$.

Proof: In the free group $F(\mathcal{R}[t, t^{-1}])$, set $x = y = 1$ (i.e. factor by the ideal $\Sigma$ in $\mathcal{R}$). What remains is an infinite cyclic group (generated by the matrix we called T at the beginning of this section). Hence, the commutator subgroup of F is in the kernel of this map and the result follows.

Remarks:
(1) This proof generalizes readily when the rank is greater than two. (For rank r, setting $x_1 = \ldots = x_r = 1$ produces a free abelian group of rank r-1. We will return to this in a later paper.)

(2) If f is in F'' the second derived subgroup of $F(\mathcal{R}[t, t^{-1}])$, then one can show that $A_1$, the coefficient of t-1, is in $\Sigma^3$. However the remaining coefficients are in $\Sigma^2$. For the next lemma, we will be content to merely assume that the coefficients of f in F'' are in $\Sigma$ as asserted in Lemma 5. Using lemma 5 as a starting point, we show that as we move down the derived series of $F(\mathcal{R}[t, t^{-1}])$, we double the power of the augmentation ideal in the entries of the coefficients of the $(t-1)^i$. The proof of Lemma 6(ii) follows that of Vaughan-Lee.

Lemma 6: (i) g is in F'' if and only if g has the form
$$g = I + (t-1)A_1 + (t-1)^2 A_2 + (t-1)^3 A_3 + \cdots$$

(ii) If g is in $F^{(k)}$, the $k^{th}$ derived group of F, (k > 1), and $d = 2^{k-2}$, then
$$g = I + (t-1)^d A_d + (t-1)^{d+1} A_{d+1} + \cdots$$
where $A_d, A_{d+1}, A_{d+2}, \ldots$ are in $\Sigma^d$.

Proof: (i) is a restatement of Lemma 1(i).

(ii) The proof is by induction on k. The start of our induction is k = 2 which is covered by Lemma 5 and part (i) above. Suppose therefore that the result is true for $k \geq 2$ and let $d = 2^{k-2}$. Let g, h be in $F^{(k)}$

$$g = I + (t-1)^d A_d + (t-1)^{d+1} A_{d+1} + \ldots \quad , \quad h = I + (t-1)^d B_d + (t-1)^{d+1} B_{d+1} + \ldots \quad ,$$

and consider [g, h]. By induction we assume that $A_i$, $B_i$ are in $\sum^d$ for $i \geq d$. Since $g^{-1}$, $h^{-1}$ are also in $F^{(k)}$, we we can write

$$g^{-1} = I + (t-1)^d C_d + (t-1)^{d+1} C_{d+1} + \ldots \quad , \quad h^{-1} = I + (t-1)^d D_d + (t-1)^{d+1} D_{d+1} + \ldots \quad ,$$

and $C_i$, $D_i$ are in $\sum^d$ for $i \geq d$. With the summations over $i \geq d$, let

$$A = \sum (t-1)^i A_i \, , \quad B = \sum (t-1)^i B_i \, , \quad C = \sum (t-1)^i C_i \, , \quad D = \sum (t-1)^i D_i \, .$$

Then $g\, g^{-1} = (I + A)(I + C) = I + A + C + AC$, and hence $A + C + AC = 0$. Similarly $B + D + BD = 0$. Hence

$$[g, h] = (I + A)(I + B)(I + C)(I + D) =$$
$$I + A + B + C + D + AB + AC + AD + BC + BD + CD + ABC + ABD + ACD + BCD + ABCD =$$
$$I + AB + AD + BC + CD + ABC + ABD + ACD + BCD + ABCD.$$

So if we expand [g, h] as a power series

$$[g,h] = I + (t-1)E_1 + (t-1)^2 E_2 + (t-1)^3 E_3 + \cdots ,$$

then $E_i = 0$ for $i < 2d$, and for $i \geq 2d$, $E_i$ is a linear combination of products of two or more matrices from the set $\{A_d, A_{d+1}, \ldots, B_d, B_{d+1}, \ldots, C_d, C_{d+1}, \ldots, D_d, D_{d+1}, \ldots\}$. By induction, all matrices in this set are in $\sum^d$, hence the $E_i$ are in $\sum^{2d}$ for $i \geq 2d$. This completes the proof.

Lemma 6(ii) is enough for our purposes as it will enable us to show that the group $F(S[t, t^{-1}])$ defined below is solvable. However, this gives a bound for the derived length of $F(S[t, t^{-1}])$ one larger than best possible since we assumed that an element of F" has all coefficients of $(t-1)^i$ in $\sum$, when it is easy to see that they lie in $\sum^2$. For completeness, we will use Lemma 7 which will give the best possible bound.

Lemma 7: If g is in $F^{(k)}$, the $k^{th}$ derived group of F, (k > 1), and $d = 2^{k-2}$, then

$$g = I + (t-1)^d A_d + (t-1)^{d+1} A_{d+1} + \cdots$$

where $A_d, A_{d+1}, \ldots$ are in $\sum^{2d}$.

Proof: Because of Lemma 6(ii), we only have to show that an element in F" has all coefficients in $\sum^2$. To show this, we use Lemma 5 as a starting point and proceed as in the proof of Lemma 6(ii). The only difference is that the constant term of an element in F', instead of the identity matrix, now has the form $u + [\lambda_i v]$, where $\lambda_1$ and $\lambda_2$ are in $\sum$. Otherwise the proof proceeds exactly as in Lemma 6(ii). We omit the details.

Recall from the previous section that $S = S(n) = \mathcal{R}/\mathcal{I}(n)\sum$. For the rest of this section we require that $q = p^e$ is a prime power so that we may apply Lemma 3(i) of the previous section. The proof of item (**iv**) of the introduction is contained in Lemma 8.

Lemma 8: (i) If $n = q = p^e$ is a prime power, $F(S[t, t^{-1}])$ is solvable.

    (ii) The derived length of $F(S[t, t^{-1}])$ is at most k where $2^{k-1} \geq e(p^e - p^{e-1}) + 1$.

Proof: We apply Lemma 7 to $F(\mathcal{R}[t, t^{-1}])$ before mapping down to $F(S[t, t^{-1}])$. Choose k large enough so that $2^{k-1} \geq e\phi(q) + 1 = e(p^e - p^{e-1}) + 1$. Then all the coefficients of $(t-1)^i$ are the zero matrix in an expansion of an element of the kth derived group of $F(S[t, t^{-1}])$. Thus elements in the $k^{th}$ derived group are the identity. This completes the proof of Lemma 8.

Since $\Sigma^{e\phi(q)-1} \not\subset \mathcal{I}(q)$, it is easy to see that the value of k in Lemma 8(ii) is best possible. Although $F(S[t, t^{-1}])$ is solvable, a difficult exercise will show that it is not nilpotent. (Factoring by $(t-1)^i$ for suitable i maps $F(S[t, t^{-1}])$ onto a nilpotent group.)

## 4.     The group $G = G(q) = F(S[t, t^{-1}])$.

Assume that $S = S(q)$, where $q = p^e$ is a prime power, and hence $G = F(S[t, t^{-1}])$ is a solvable group. (The solvability class depends on q and tends to infinity as q does.) The generators of G are the images of $M_1$ and $M_2T$ under $\alpha$. To simplify notation we will use the same letters for the elements of G as for $F(\mathcal{R}[t, t^{-1}])$ but with a line through them. Hence the lined letters are matrices over $S[t, t^{-1}]$. Thus
$$G = F(S[t, t^{-1}]) = gp \langle \alpha(M_1), \alpha(M_2T) \rangle = gp \langle \overline{M}_1, \overline{M}_2T \rangle$$
( T is left unchanged since $\alpha$ leaves T unchanged.)

Clearly $\overline{M}_1$ has order q and it is easy to see that $\overline{M}_2T$ has infinite order. Namely, sending x and y to 1 sends $\overline{M}_2T$ to T which is a matrix of infinite order. (Without danger of confusion we are using x and y as elements of $S$ as well as $\mathcal{R}$.) If $\overline{W}$ is an element of G, i.e., a product of $\overline{M}_1$, $\overline{M}_2T$ and their inverses, such that $\overline{W}$ has a nonzero exponent sum in T, then the same proof shows that $\overline{W}$ has infinite order. Next consider the elements with zero exponent sum in T. Examples are any conjugate of $\overline{M}_1$ which of course have order q. However, there is a much better way to show this, relying heavily on Theorem B in [3] ( i.e., Lemma 3 (i)). Such a proof shows that any element which has exponent sum zero in T has order a divisor of q. We state this as Theorem 2 and limit our proof to primes, an enormous simplification over prime powers.

Theorem 2: Let $\overline{W}$ be an element of $G = F(S[t, t^{-1}])$.
a) If $\overline{W}$ has exponent sum zero in T, then $\overline{W}^q = 1$.
b) If $\overline{W}$ has a nonzero exponent sum in T, then $\overline{W}$ has infinite order.

Proof: Assume q is a prime. If $\overline{W}$ has exponent sum zero in T, then we claim that $\overline{W}^q = 1$. $\overline{W}$ is a product of $\overline{M}_1$, $\overline{M}_2$T and their inverses. Let W be a preimage of $\overline{W}$ in $F(\mathcal{R}[t, t^{-1}])$. Write the $M_i$ in the form of Lemma 1(ii) of section 2 and write T in the form $tI + A$, where tI is the scalar matrix with t in the diagonal. As a product of these matrices W has the form $W = uI + V$, where $u = x^i y^j t^k$ and the entries of V are in the augmentation ideal of $\mathcal{R}$. By our assumption, $k = 0$ (the only place we use the hypothesis that the exponent sum of T is zero). Now map W to $\overline{W}$, i.e., send $\mathcal{R}$ to $S$. Since $S$ has prime characteristic q, $(uI + \overline{V})^q = u^q I + \overline{V}^q$, where $u^q = 1$ and $\overline{V}^q$ is the zero matrix. Namely, the entries of $\overline{V}^q$ are in the qth power of the augmentation ideal. By Lemma 3(i), the q-1 power of the augmentation ideal lies in $\mathcal{I}(q)$ and hence the qth power is in $\mathcal{I}(q)\Sigma$ which is the zero of S.

Combining Theorem 2 with Theorem 1 presents a more complete picture. Recall that setting $t = 1$ in a matrix of G sends that matrix into one whose order is a divisor of q. (Factoring by the ideal $t - 1$ in $S[t, t^{-1}]$ is factoring by the second derived group in G.) Thus as a corollary of Theorems 1 and 2, we have the

Theorem 3: Let g be an element in $G = G(q) = F(S[t, t^{-1}])$. Then either $g^q = 1$ or g has infinite order. If g has infinite order, $g^q \equiv 1$ modulo G''.

The ring $S = S(q) = \mathcal{R}/\mathcal{I}(q)\Sigma$ was chosen as precisely the ring which changes the free metabelian group $F(\mathcal{R})$ into $F(S)$, the Burnside metabelian group of exponent q. Using the same ring map to change $F(\mathcal{R}[t, t^{-1}])$ into $F(S[t, t^{-1}])$ creates a substantial number of elements of exponent q in $F(S[t, t^{-1}])$, although the extent as described in Theorem 2 is perhaps surprising. But it does appear to give a good group theoretic answer to the question as to *why should $F(S[t, t^{-1}])$ be solvable*. Namely, enough elements of exponent q imply (higher) commutator identities. Unfortunately no one, to my knowledge, has any idea what these identities look like when q is 5 or larger.

## 5. A further application of Theorem 1

Consider the Burnside map $\beta: F(\mathcal{R}[t, t^{-1}]) \to B(q)$, where B(q) is a group isomorphic to $F/F^q$. Since $F(\mathcal{R}[t, t^{-1}]$ is a free group, we know such a map exists even if we are not sure what form B(q) might take. (As will become clear shortly, there is overwhelming evidence that B(q) is in fact a matrix group over a Noetherian ring.) It is however trivial to observe that the generator $M_1$ of $F(\mathcal{R}[t, t^{-1}])$, which does not contain t in any entry, is mapped by β into a matrix over $S$ (see lemma 10 below). Our motivation is to show that B(q) is a homomorphic image of $F(S[t, t^{-1}]$ and hence a solvable group. Notice that Theorem 2 already suggests this, but we will not pursue a proof in this direction and will not use Theorems 2 or 3.

For convenience, we begin by displaying the diagram (pg. 1) which contains the groups $F(\mathcal{R}[t, t^{-1}])$, $F(\mathcal{R})$, $F(S[t, t^{-1}])$ and $F(S)$ and their maps $\alpha$ and $\underline{\alpha}$.

$$\begin{array}{ccc} F(\mathcal{R}[t, t^{-1}]) & \xrightarrow{\alpha} & F(S[t, t^{-1}]) \\ \downarrow & & \downarrow \\ F(\mathcal{R}) & \xrightarrow{\underline{\alpha}} & F(S) \end{array}$$

For g in $F(\mathcal{R}[t, t^{-1}])$, we denote by g/t the matrix obtained by setting t = 1 in g. If $\lambda$ is a map of $F(\mathcal{R}[t, t^{-1}])$ onto a group G, then g/t is an element of $F(\mathcal{R})$ and thus, in general, is not in the domain of $\lambda$. But one can define a new map $\underline{\lambda}$ on $F(\mathcal{R})$ if one can unambiguously evaluate $\lambda(g/t)$. We can do this because $F(\mathcal{R})$ is free metabelian. (One can define a valid map on any quotient group of the free group $F(\mathcal{R}[t, t^{-1}])$ by a fully invariant subgroup.) Thus it is easy to show:

Lemma 9: $\underline{\lambda}: F(\mathcal{R}) \to G/G''$ given by $\underline{\lambda}(g/t) = \lambda(g)G''$ is a well defined map of $F(\mathcal{R})$ onto $G/G''$. $\underline{\lambda}$ is the map *induced* from $\lambda$.

Proof: If it is well-defined, then clearly $\underline{\lambda}$ is a homomorphism. Suppose therefore that g/t and h/t are equal elements of $F(\mathcal{R})$. Then $h^{-1}g/t$ is the identity and thus $h^{-1}g$ is in the second derived group of $F(\mathcal{R}[t, t^{-1}])$. But then $\lambda(h^{-1}g)$ is in $G''$ since $\lambda$ takes the second derived group of $F(\mathcal{R}[t, t^{-1}])$ into $G''$. Thus $\underline{\lambda}(g/t) = \lambda(g)G'' = \lambda(h)G'' = \underline{\lambda}(h/t)$.

example: In Lemma 9, suppose that $G = F(S[t, t^{-1}])$ and $\lambda = \alpha$. The kernel of $F(S[t, t^{-1}]) \to F(S)$, which is our $G''$, consists of the matrices whose entries are in the ideal in $S[t, t^{-1}]$ generated by t-1. Again we need to compute when t = 1, but now in the image of $\alpha$. Thus the *induced* map $\underline{\alpha}: F(\mathcal{R}) \to F(S)$ is given by

$$\underline{\alpha}(g/t) = \alpha(g)/t \quad \text{for g in } F(\mathcal{R}[t, t^{-1}]).$$

Notice that although $\underline{\alpha}$ was defined by the ring homomorphism $\mathcal{R} \to S$, we could have defined $\underline{\alpha}$ as the map on $F(\mathcal{R})$ *induced* by $\alpha$. Once one knows a map which induces a given map, it is usually easy to find others that induce the same map. In addition to $\alpha$, there are an infinite number of maps on $F(\mathcal{R}[t, t^{-1}])$ which induce $\underline{\alpha}$. But all will have to contain enough qth powers in their image since $\underline{\alpha}$ is the Burnside map on the free metabelian group.

Remarks:

(1) For abstract groups, the usual definition of induced map is as follows:
   If F is a free group and $\lambda: F \to G$ is a map, the map $\underline{\lambda}: F/F'' \to G/G''$ is defined by $\underline{\lambda}(gF'') = \lambda(g)G''$ for g in F. One says that $\lambda$ *induces* $\underline{\lambda}$ or that $\underline{\lambda}$ is *induced* by $\lambda$. Rather than an abstract free group, we started with a free group of matrices $F(\mathcal{R}[t, t^{-1}])$. Our definition was adapted accordingly. As done in an earlier draft, one could begin with the abstract definition and derive our definition as a consequence.

(2) The maps on $F(\mathcal{R}[t, t^{-1}])$ that concern us are $\alpha$ and the Burnside map $\beta$. Both maps *induce* $\underline{\alpha}$. Namely, $\beta: F(\mathcal{R}[t, t^{-1}]) \to B(q)$ *induces* $\underline{\beta}$ which maps $F(\mathcal{R})$ onto $B(q)/B(q)''$, which is isomorphic to the metabelian Burnside group $F(S)$.

In what follows, bear in mind that $\underline{\alpha}$ sends the ring $\mathcal{R}$ into the ring $S$. We will use the same letters for the elements of $F(S[t, t^{-1}])$ or $F(S)$ as those used for $F(\mathcal{R}[t, t^{-1}])$ or $F(\mathcal{R})$ but with a line through them. Thus, lined letters are always matrices over $S$, although $S$ may be embedded in $S[t, t^{-1}]$.

Lemma 10: Suppose $\lambda$ induces $\underline{\alpha}$. Then $\lambda(M_1) = \alpha(M_1) = \overline{M}_1$.
Proof: By assumption, $\underline{\lambda}(M_1/t) = \underline{\alpha}(M_1/t)$. Since $M_1$ does not contain t in any entry, $M_1/t = M_1$. Therefore $\underline{\lambda}(M_1/t) = \underline{\lambda}(M_1)$ and $\underline{\alpha}(M_1/t) = \underline{\alpha}(M_1) = \overline{M}_1$. Hence $\lambda(M_1)$ and $\alpha(M_1)$ are the same, both being equal to $\overline{M}_1$.

Remark: If g is a matrix of $F(\mathcal{R}[t, t^{-1}])$ without t in any entry, then Lemma 10 merely states that a map $\lambda$ on $F(\mathcal{R}[t, t^{-1}])$ has the same image $\lambda(g)$ as the "specialization" $\lambda(g/t)$ at t = 1, which is the induced map $\underline{\lambda}$.

Lemma 11: Suppose $\lambda$ induces $\underline{\alpha}$ and g is an element of $F(\mathcal{R}[t, t^{-1}])$. Then in the entries of g, $\lambda$ sends the elements of $\mathcal{R}$ into elements of $S$.
Proof: By Lemma 10, $\lambda$ sends the entries of $M_1$ into the ring $S$. Since $M_1$ contains the generators of $\mathcal{R}$, $\lambda$ sends the elements of $\mathcal{R}$ into $S$ in all the matrices of $F(\mathcal{R}[t, t^{-1}])$.

Taking $\beta$ for $\lambda$, we restate Lemma 11 as

Theorem 4: Suppose g is an element of $F(\mathcal{R}[t, t^{-1}])$. Then in the entries of g, $\alpha$ and $\beta$ are the same on $\mathcal{R}$.; i.e., both send the elements of $\mathcal{R}$ into their respective elements of $S$.

Remarks:
1) Suppose g is an element of $F(\mathcal{R}[t, t^{-1}])$. Then g is a product of $M_1$, $M_2T$ and their

inverses. Since all entries of the $M_i$ are in $\mathcal{R}$, $\alpha$ and $\beta$ differ on $F(\mathcal{R}[t, t^{-1}])$ only by their action on t. Since $\alpha$ fixes the entries of T (and hence fixes T ), the relations in the image of $\alpha$ (i.e., in $F(S[t, t^{-1}])$) are determined by the change of ring from $\mathcal{R}$ to $S$. Theorem 4 says that the image of $\beta$ contains at least those relations (and surely further ones since T is not fixed by $\beta$). Thus $F/F^q$ is a solvable group.

2)   For the generator $M_1$ of $F(\mathcal{R}[t, t^{-1}])$, $\beta(M_1)$ is the matrix $\tilde{M}_1$ of order q with entries in $S$. This, and the first remark, is a rather overwhelming indication that the image of $\beta$ is a matrix group over a quotient of $\mathcal{R}[t, t^{-1}]$, i.e., a Noetherian ring ( in fact a quotient of $S[t, t^{-1}]$). If this is the case, then results in Wehrfritz [9] chapter 11 enable one to conclude that the image of $\beta$ is a finite group.

## Appendix

Lemma 1(ii):  Each element M of $F(\mathcal{R})$ has the form  $M = uI + N$, where I is the identity matrix, u is a positive unit in $\mathcal{R}$, and $N = [\lambda_i v]$ is an r-square matrix whose $i^{th}$ row is $\lambda_i v$. The $\lambda_i$ in $\mathcal{R}$ satisfy $\lambda_1(1-x_1) + \cdots + \lambda_r(1-x_r) = 1- u$.

Proof:  It is a simple matter to check that $vM = v$ for all M in $F(\mathcal{R})$. (i.e. verify this for each generator in $F(\mathcal{R})$.)  Also verify that each generator and it's inverse satisfies Lemma 9(ii). Having verified Lemma 9(ii) for words of length 1, we use induction on the length of a word in $F(\mathcal{R})$. Suppose $W_1W_2$ is an element of $F(\mathcal{R})$, where

$W_1 = u_1+ [\lambda_i v]$, $W_2 = u_2+ [\delta_i v]$ and $\lambda_1(1-x_1) + \cdots + \lambda_r(1-x_r) = 1- u_1$,

$\delta_1(1-x_1) + \cdots + \delta_r(1-x_r) = 1- u_2$. Then using the fact that $vM = v$ for all M in $F(\mathcal{R})$, $W_1W_2 = (u_1+ [\lambda_i v])(u_2+ [\delta_i v]) = u_1(u_2+ [\delta_i v]) + [\lambda_i v] = u_1u_2+ u_1[\delta_i v]) + [\lambda_i v]$.

But $(u_1\delta_1+ \lambda_1)(1-x_1) + \cdots + (u_1\delta_r + \lambda_r)(1-x_r) = u_1(1- u_2)+ (1- u_1)= 1- u_1u_2$. Thus the result holds for $W_1W_2$